\documentstyle[11pt]{article}
\newtheorem {th}{Theorem}
\newtheorem {lem}[th]{Lemma}

\def\Cox{\hfill \Box}

\def\sf{\sigma\mbox{-field}}

\def\ee{\epsilon}
\def\E{{\bf{E}}}
\def\P{{\bf{P}}}
\def\N{\hbox{I\kern-.2em\hbox{N}}}
\def\R{\hbox{I\kern-.2em\hbox{R}}}
\def\Z{{\bf{Z}}}

\def\C{{\cal{C}}}

\def\F{{\cal{F}}}
\def\|{\, | \, }

\def\vv{{\bf v}}
\def\vn{{\bf v^{(n)}}}

\begin{document}

\begin{center}
{\large \bf Planar first-passage percolation times are not tight}
\end{center}

\begin{flushright}
R. Pemantle\footnote{Department of Mathematics, University of 
Wisconsin-Madison, Van Vleck Hall, 480 Lincoln Drive, Madison, 
WI 53706}$^,$\footnote{This research supported in part by NSF grant 
\#DMS-9300191, by a Sloan Foundation Fellowship and by a 
Presidential Faculty Fellowship.} and Y. Peres\footnote{Department
of Statistics, University of California at Berkeley, 367 Evans Hall,
Berkeley, CA 94720}
\end{flushright}

We consider first-passage percolation on the two-dimensional integer
lattice $\Z^2$ with passage times that are IID exponentials of mean one;
see Kesten (1986) for an overview. 
It has been conjectured, based on numerical evidence, that the
variance of the time $T(0,n)$ to reach the vertex $(0,n)$ is of
order $n^{2/3}$.  Kesten (1992) showed that the variance of $T(0,n)$ is 
at $O(n)$.  He also noted that the variance is
bounded away from zero.  This note improves the lower bound on the variance
of $T(0,n)$ to $C \log n$.  Simultaneously and independently,
Newman and Piza have achieved the same result for 
$\{ 0 , 1 \}$-valued passage times.  Their methods (Newman and Piza 1993) 
extend to more general passage times, while ours work only for exponential
times.  On the other hand, our theorem shows that the variance 
comes from fluctuations of nonvanishing probability in the sense 
that, as $n \rightarrow \infty$, the law of $T(0,n)$ is not tight 
about its median.  Very recently, Newman and Piza showed that
the $\log n$ may be improved to a power of $n$ 
for directions in which the shape is not flat (it is not known whether
the shape can be flat in any direction; see Durrett and Liggett
(1981) for a relevant example).  As pointed out to us by Harry Kesten, 
in the exponential case this may also be obtained via the method
given here.
 
\begin{th}
Let $\vv$ be any unit vector in $\R^2$ and let $\vn$ be the vector
in $\Z^2$ whose coordinates are the integer parts of the
coordinates of $n\vv$.  Let $T(0,\vn)$ denote the passage
time from the origin to the vertex $\vn$, under IID 
mean 1 exponential passage times on the edges of $\Z^2$.  Then
$\mbox{Var}(T(0,\vn)) \geq C \log n$ and in fact
any intervals $[a_n , b_n]$ with $b_n - a_n = o(\log n)^{1/2}$ satisfy
$\P (T(0,\vn) \in [a_n , b_n]) \rightarrow 0$.
\end{th}

\noindent{Remark:} The theorem extends to Richardson's model with
other passage time distributions (restart all clocks after an edge is
crossed).

\noindent{\sc Proof:}  We compute the conditional distribution
of $T(0,\vn)$ given a $\sf$ $\F$ and show that with probability
$1 - o(1)$, this conditional distribution is close to a normal
with variance at least $C \log n$; clearly this implies
the conclusion of the theorem.

Let $\F$ be the $\sf$ determined by the order in which vertices are
reached.  Formally, if $T(v)$ is the passage time from the origin to
the vertex $v$, then $\F$ is the $\sf$ generated by the events 
$T(v) < T(w)$ for $v,w \in \Z^2$.  Let $V_0 , V_1 , V_2 , \ldots $ be the
vertices of $\Z^2$ listed in the order they are reached, so 
$(V_0 , V_1 , \ldots)$ is an $\F$-measurable random sequence.
Let $\C_n = \{ V_0 , V_1 , \ldots , V_n \}$ be the cluster of the first
$n$ elements to be reached from the origin, and let $Y_n$ be the
number of edges connecting elements $\C_{n-1}$ to elements of $\Z^2
\setminus \C_{n-1}$.  The key observation is that the conditional
joint distribution of the variables $T(V_n) - T(V_{n-1})$ given
$\F$ is identical to a sequence of independent exponentials with
means $1 / Y_n$.  This is in fact an immediate consequence of the
lack of memory of the exponential distribution and of the fact that
the minimum of $n$ independent exponentials of mean 1 is an exponential
of mean $1/n$.  This observation leads to
\begin{lem} \label{lem1}
$$ \liminf_{n \rightarrow \infty} (\log n)^{-1} \sum_{j=1}^n 
   {1 \over Y_j^2} \geq c_0 ~\mbox{ a.s.}$$
for some positive constant $c_0$.
\end{lem}
Assuming this lemma for the moment, define $M(n)$ to satisfy 
$V_{M(n)} = \vn$.  Let 
$$\mu_n = \E (T(V_n) \| \F) = \sum_{j=1}^n {1 \over Y_j} $$ 
and 
$$\sigma_n^2 = \mbox{Var} (T(V_n) \| \F) = \sum_{j=1}^n {1 \over Y_j^2} .$$
The Lindeberg-Feller theorem implies that the conditional 
distribution of the variable \\
$(T(\vn) - \mu_{M(n)}) / \sigma_{M(n)}$ converges weakly to a 
standard normal whenever $\sigma_{M(n)} \rightarrow \infty$.  
Subadditivity implies that
that $T(\vn) / n \rightarrow c_1 = c_1 (\vv)$ almost surely, and the
shape theorem (Cox and Durrett 1981) implies that $c_1 > 0$ and
that the number of vertices $N_t$ reached by time $t$ is almost surely
$(c_2 + o(1)) t^2$.  (Equivalently, $T(V_n) = (c_2^{-1/2} + o(1)) 
n^{1/2}$.)  From this it follows that
\begin{equation} \label{eq1}
{M(n) \over n^2} \rightarrow c_2 c_1^2 \mbox{~~a.s.},
\end{equation}
and hence from Lemma~\ref{lem1} that 
$$\liminf (\log n)^{-1} \sigma_{M(n)}^2 \geq 2 c_0 \mbox{~~a.s.}$$
Thus a conditional distribution of $T(\vn)$ is close to a normal
with variance at least $2 c_0 + o(1)$, which establishes the
theorem.

To prove Lemma~\ref{lem1}, we first observe from the isoperimetric
inequality that the distribution of $T(V_n)$ is stochastically
bounded by the sum of independent exponentials of means 
$c j^{-1/2}$, $j = 1 , \ldots , n$.  Thus the variables $n^{-1/2} T(V_n)$
are dominated by a variable in $L^1$ and hence
$${\mu_n \over n^{1/2}} = \E (n^{-1/2} T(V_n) \| \F) \rightarrow 
   c_2^{-1/2}$$
almost surely.  Lemma~\ref{lem1} then follows from a fact about
sequences of real numbers:

\begin{lem} \label{lem2}
Let $x_1 , x_2 , \ldots$ be positive real numbers with $S_n = 
\sum_{j=1}^n x_j$ and suppose that \\ $\liminf n^{-1/2} S_n = c$.  Then 
$$\liminf (\log n)^{-1} \sum_{j=1}^n x_j^2 \geq {c^2 \over 4} .$$
\end{lem}

\noindent{\sc Proof:}  It suffices to show that the condition
$$S_n \geq a n^{1/2} - b ~\mbox{ for all } n$$
implies $\liminf (\log n)^{-1} \sum_{j=1}^n x_j^2 \geq a^2 / 4$, since
one may then take $a = c - \ee$ for arbitrarily small $\ee$.  Also,
replacing $x_1$ with $x_1 + b$, we may assume without loss of generality
that $b = 0$.  Define
$$q_n = n^{1/2} - (n-1)^{1/2} \geq {1 \over 2} n^{-1/2} .$$
Assuming $S_n \geq a n^{1/2}$ for all $n$, we show that
$\sum_{j=1}^n x_j^2 \geq a^2 \log n / 4$.  
Rearranging the terms $\{x_j : 1 \leq j \leq n \}$ in decreasing
order does not change $\sum_{j=1}^n x_j^2$ and only increases each $S_j$,
so we may assume without loss of generality that these terms appear in
decreasing order.  Summing by parts three times we obtain
\begin{eqnarray*}
\sum_{j=1}^n x_j^2 & = & S_n x_n + \sum_{k=1}^{n-1} S_k 
   (x_k - x_{k+1}) \\[2ex]
& \geq & a \left [ n^{1/2} x_n + \sum_{k=1}^{n-1}k^{1/2} (x_k - x_{k+1})
   \right ] \\[2ex]
& = & a \sum_{j=1}^n q_j x_j \\[2ex]
& = & a \left [ q_n S_n + \sum_{k=1}^{n-1} (q_k - q_{k+1}) S_k \right ]
   \\[2ex]
& \geq & a \left [ q_n n^{1/2} + \sum_{k=1}^{n-1} (q_k - q_{k+1})
   k^{1/2} \right ] .
\end{eqnarray*}
Summing once more by parts and using the definition of $q_k$ we
see that this is equal to $a \sum_{k=1}^n q_k^2$, which is
at least $a^2 \log n / 4$.  This proves the lemma and hence 
the theorem.  $\Cox$

When the asymptotic shape has a finite radius of curvature in 
the direction $\vv$, Newman and Piza have shown, using results
of Kesten (1992) and Alexander (1992), that the minimizing path from 
$0$ to $\vn$ deviates from a straight line segment by at most
$c n^\alpha$ for some $\alpha < 1$, with probability $1 - o(1)$
as $n \rightarrow \infty$.  Thus the time $T'(0,\vn)$ to reach
$\vn$, in a new percolation where only bonds in a strip of width 
$n^\alpha$ are permitted, differs from $T(0,\vn)$ by $o(1)$ in total
variation.  The shape theorem for $\Z^2$ implies that the number
of sites reached in the new percolation by this time is 
$O(n^{1+\alpha})$.  Defining $M'(n), \mu_n', Y_n'$ and $\sigma_n'$
analogously to $M(n), \mu_n , Y_n$ and $\sigma_n$
but for the new percolation, we have $\mu_n' = \sum_{k=1}^n (1/ Y_n')$,
while now $M'(n) = O(n^{1+\alpha})$.  By Cauchy-Schwartz
-- no summing by parts is needed -- it follows that
$(\sigma_{M'(n)}')^2 \geq c n^{1-\alpha}$, and applying Lindeberg-Feller
as before proves the extension mentioned before Theorem~1.


\begin{thebibliography}{YMN}

\bibitem{AK}
Alexander, K. (1992).  Fluctuations in the boundary of the wet region
for first-passage percolation in two and three dimensions.  {\em Preprint.}

\bibitem{CD}
Cox, J.T. and Durrett, R. (1981).  Some limit theorems for percolation 
processes with necessary and sufficient conditions.  {\em Ann. Probab.}
{\bf 9} 583 - 603.

\bibitem{DL}
Durrett, R. and Liggett, T. (1981).  The shape of the limit set in 
Richardson's growth model.  {\em Ann. Probab.} {\bf 9} 186 - 193.

\bibitem{Ke1}
Kesten, H. (1986).  Aspects of first passage percolation.  In Lecture
Notes in Mathematics, vol. 1180, 125 - 264.

\bibitem{Ke2}
Kesten, H. (1992).  On the speed of convergence in first passage 
percolation.  {\em Ann. Appl. Prob. to appear}

\bibitem{NP}
Newman, C. and Piza, M. (1993).  Divergence of shape fluctuations
in two dimensions.  {\em Preprint.}

\end{thebibliography}
\end{document}